\documentclass[aps,prl,reprint,groupedaddress]{revtex4-2}
\usepackage{amsmath}
\usepackage{bm}
\usepackage{amsfonts}
\usepackage{latexsym}
\usepackage{amssymb}
\usepackage{stmaryrd}
\usepackage{overpic}
\usepackage{graphicx}
\usepackage{multirow}
\usepackage{mathrsfs}
\usepackage{float}
\usepackage{color}
\usepackage{setspace}
\usepackage{subfigure}
\usepackage{dsfont}

\begin{document}

\title{FEM-DtN-SIM Method for Computing Resonances of Schr\"{o}dinger Operators}
\author{Bo Gong$^1$, Takumi Sato$^2$, Jiguang Sun$^3$, Xinming Wu$^2$}
\affiliation{$^1$School of Mathematics, Statistics and Mechanics, Beijing University of Technology, Beijing, 100124, China.\\$^2$School of Mathematical Sciences, SKLCAM, Fudan University, Shanghai 200433, China. \\$^3$Department of Mathematical Sciences, Michigan Technological University, Houghton, MI 49931, U.S.A.}

\begin{abstract}
The study of resonances of the Schr\"{o}dinger operator has a long-standing tradition in mathematical physics. Extensive theoretical investigations have explored the proximity of resonances to the real axis, their distribution, and bounds on the counting functions. However, computational results beyond one dimension remain scarce due to the nonlinearity of the problem and the unbounded nature of the domain. We propose a novel approach that integrates finite elements, Dirichlet-to-Neumann (DtN) mapping, and the spectral indicator method. The DtN mapping, imposed on the boundary of a truncated computational domain, enforces the outgoing condition. Finite elements allow for the efficient handling of complicated potential functions. The spectral indicator method effectively computes (complex) eigenvalues of the resulting nonlinear algebraic system without introducing spectral pollution. The viability of this approach is demonstrated through a range of numerical examples.
\end{abstract}

\maketitle
\section{1. Introduction}
Resonances of Schr\"{o}dinger operators have been a classical topic in mathematics and physics, generalizing bound states for systems where energy can scatter to infinity. They are characterized as poles of the meromorphic continuation of the resolvent operator, with their real parts representing oscillation frequencies and their imaginary parts corresponding to decay rates \cite{DyatlovZworski2019}. Extensive studies have investigated the proximity of resonances to the real axis, resonance-free regions, bounds, and counting functions, etc. \cite{BrietCombesDuclos1987JMAA, DyatlovZworski2019, LiYaoZhao2024SIMA}. The study of scattering resonances has important applications in various fields, including gravitational wave signatures, photonics, and plasmonic scattering \cite{Maakitalo2014, JMS2022PRL, BBCGB2024PRB}.

In contrast to the extensive literature on the analysis of resonances, relatively few studies address the computation of resonances of Schrödinger operators beyond the one-dimensional case \cite{BindelZwordki2006, BMR2023}. The problem becomes significantly more challenging in higher dimensions. When restricting the computation to a bounded domain, it is crucial to impose the outgoing condition correctly. Additionally, the problem exhibits a nonlinear dependence on the eigen-parameter. In one dimension, the outgoing condition results in a quadratic eigenvalue problem, which can be linearized and efficiently solved using standard numerical linear algebra techniques \cite{BindelZwordki2006}. This linearization is no longer feasible in higher dimensions. Moreover, formulations that employ absorbing boundary conditions, such as the perfectly matched layer (PML) or impedance boundary conditions, are not equivalent to the original problem. As a result, they can introduce significant spectral pollution or even yield entirely incorrect resonance computations.

In this paper, we study the computation of resonances of the Schr\"{o}dinger operator with a bounded (complex) potential of compact support. By reformulating the problem on a truncated domain and incorporating the Dirichlet-to-Neumann (DtN) mapping to enforce the outgoing condition, we obtain an eigenvalue problem for a holomorphic operator-valued function. We then propose a finite element method to solve this nonlinear system, where the computed eigenvalues correspond to the resonances \cite{Hsiao2011, SunZhou2017}. The discrete eigenvalue problem is computed using a parallel multistep spectral indicator method based on contour integrals \cite{XiSun2023}. We provide numerical examples for various potentials to validate the proposed method. While our focus is on two-dimensional cases, the scheme is also applicable to three-dimensional problems. The computational framework developed here serves as a quantitative tool for studying Schr\"{o}dinger resonances. Notably, it has also been successfully applied to compute scattering resonances for acoustic obstacles and metallic grating structures \cite{XiGongSun2024,XiLinSun2024CMA}.

The convergence of the proposed method can be established using abstract approximation theory for eigenvalue problems of holomorphic operator functions, following the same approach as in the analysis of scattering resonances for acoustic obstacles \cite{karma1996a, XiGongSun2024}. Specifically, the computed resonances converge to the true resonances at a certain rate as the mesh size $h \to 0$, which is further demonstrated through numerical experiments.

The remainder of the paper is organized as follows. In Section~2, we introduce the concept of resonances of the Schr\"{o}dinger operator. In Section~3, we truncate the infinite domain using a circle and impose the Dirichlet-to-Neumann (DtN) mapping on the boundary. Section~4 contains the finite element formulation and the definition of a holomorphic operator-valued function whose eigenvalues correspond to the resonances. Finally, in Section~5, we provide numerical examples to validate the proposed approach.


\section{2. Schr\"{o}dinger Resonances}\label{SRDtN} 
Let $V(x) \in L^\infty(D)$ be a potential function with compact support $D \subset \mathbb R^2$.  The Schr\"{o}dinger operator is defined as
\begin{equation}\label{SchrodingerO}
H_V:= - \triangle  + V(x).
\end{equation}
The resonance problem for the Schr\"{o}dinger operator is to find $k \in \mathbb C$ and a nontrivial outgoing solution $u$ to the Schr\"{o}dinger equation
\begin{equation}\label{Schrodinger}
H_V u(x) - k^2 u(x) = 0, \quad x \in \mathbb R^2.
\end{equation}
For any $k \in \mathbb C$, an outgoing solution $u$ to the Schr\"{o}dinger equation $H_V u(x) - k^2 u(x) = f$ satisfies \cite{Zworski1999NAMS}
\begin{equation}\label{outgoingconditon}
u(x)=\sum_{n=-\infty}^{\infty}a_nH_n^{(1)}(k r)e^{in\theta},\ \ \ |x|>r_0,
\end{equation}
where $\theta=\arg(x)$, $r_0>0$ is a constant large enough, and $H_n^{(1)}(\cdot)$ is the first kind Hankel function of order $n$.

The resonances are the poles of the resolvent operator $R(k) := (-\triangle + V -k^2)^{-1}$. Note that $R(k)$ is analytic for ${\rm Im} k \gg 0$ and can be continued as a meromorphic function on $\mathbb C$. Here we are interested in the numerical computation of the resonances $k$'s in the lower half complex plane. 

\begin{figure} [ht!]
\begin{center}
\rotatebox{0}{ \scalebox{0.8} {\includegraphics*{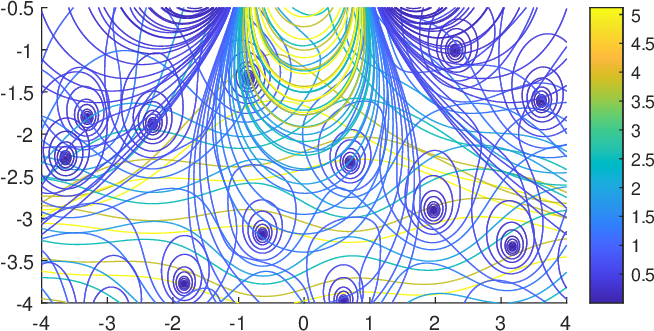}}}
\caption{Contour plots of $|d_n(k)|$, $n = 0, \ldots, 10$ ($r_0=1, V(x)=2$).}
\label{PolesSchR1V2}
\end{center}
\end{figure}

When $D$ is a disk with radius $r_0$ and $V(x)=V_0 \ge 0$, the resonances can be derived analytically. The outgoing solution $u$ can be written in the form
\[
u(x) = \left\{ \begin{array}{ll} \sum_{n=-\infty}^\infty a_n H_n^{(1)} (kr) e^{in\theta}, & r=|x|> r_0, \\ \sum_{n=-\infty}^\infty b_n J_n (\sqrt{k^2-V_0}r) e^{in\theta}, & r=|x| \le r_0. \end{array}\right.
\]
The continuity of $u$ and its derivative at the boundary of $D$ leads to 
\begin{eqnarray*}
a_n H_n^{(1)} (kr)  \Big |_{r=r_0} &=&  b_n  J_n (\sqrt{k^2-V_0}r) \Big |_{r=r_0},  \\
a_n \frac{\partial H_n^{(1)} (kr) }{\partial r} \Big |_{r=r_0}  &=&b_n \frac{\partial   J_n (\sqrt{k^2-V_0}r) }{\partial r} \Big |_{r=r_0}.
\end{eqnarray*}
The above equations can be written as a linear system
{\small
\[
\begin{pmatrix} H_n^{(1)} (kr_0) &  -J_n (\sqrt{k^2-V_0}r_0) \\   kH_{n+1}^{(1)}(kr_0)  & -\sqrt{k^2-V_0}J_{n+1}(\sqrt{k^2-V_0}r_0)\end{pmatrix} \begin{pmatrix} a_n \\ b_n \end{pmatrix} =  \begin{pmatrix} 0 \\ 0\end{pmatrix}.
\]
}
The scattering resonances are $k$'s such that
\begin{eqnarray*}
&& \sqrt{k^2-V_0}J_{n+1}\left(\sqrt{k^2-V_0}r_0\right) H_n^{(1)} (kr_0)\\&& \qquad \qquad \qquad  - k H_{n+1}^{(1)} (kr_0)  J_n \left(\sqrt{k^2-V_0}r_0\right)  =0. 
\end{eqnarray*}
If $r_0=1$, one has that
\begin{eqnarray*}
d_n(k)&:=& \sqrt{k^2-V_0}J_{n+1}\left(\sqrt{k^2-V_0}\right) H_n^{(1)} (k) \\
&& \qquad   - k H_{n+1}^{(1)} (k)  J_n \left(\sqrt{k^2-V_0}\right). 
\end{eqnarray*}
The contours of the absolute values of the above determinants $d_n(k), n=0, \ldots, 10$ are shown in Fig.~\ref{PolesSchR1V2}. The minima (zeros) indicate the locations of the exact resonances.

\section{3. Dirichlet-to-Neumann Mapping}\label{DtN} 
Since the problem is posed on the infinite domain $\mathbb R^2$, we first reformulate it as an equivalent problem on a bounded domain to facilitate finite element computation. Let $\Gamma_R$ be a circle centered at the origin with radius $R$ and $D$ inside. Denote by $\Omega$ the disk such that $\partial \Omega = \Gamma_R$. We need to impose a suitable boundary condition on $\Gamma_R$. To this end, we employ the DtN mapping on $\Gamma_R$. Let $H^1(\Omega)$ be the Sobolev space of integrable functions defined on $\Omega$ with integrable gradients.  Let $H^{1/2}(\Gamma_R)$ be the trace space of $H^1(\Omega)$ and $H^{-1/2}(\Gamma_R)$ be its dual space. Following \cite{Hsiao2011}, the DtN operator $T(k):H^{1/2}(\Gamma_R)\rightarrow H^{-1/2}(\Gamma_R)$ is defined as 
\[
    T(k)\varphi = \sum_{n=0}^{+\infty}{\!}^{{}^{\,\scriptstyle{\prime}}}\frac{k}{\pi}\frac{{H^{(1)}_n}^{\prime}(kR)}{H^{(1)}_n(kR)}\int_0^{2\pi}\varphi(\phi)\cos(n(\theta-\phi))\text{d}\phi,
\]
where the summation with a prime $'$ means that the zero-th term is factored by $1/2$. 

The operator $T(k)$ is bounded. For real number $s\geq1/2$,
$T(k):H^{s}(\Gamma_R)\rightarrow H^{s-1}(\Gamma_R)$, satisfies
\begin{equation}\label{Tk}
\|T(k)\varphi\|_{H^{s-1}(\Gamma_R)}\leq C\|\varphi\|_{H^{s}(\Gamma_R)},
\end{equation}
where $C>0$ is a constant depending on $kR$ but not $\varphi$ for almost all $k\in \mathbb C$.

Using $T(k)$, we obtain a problem defined on the bounded domain $\Omega$, which is equivalent to \eqref{Schrodinger}: Find $k \in \mathbb C$ and nontrivial $u\in H^1(\Omega)$ such that
\begin{equation}\label{SchOmega}
\left\{
\begin{array}{rlll}
-\Delta u + V(x) u - k^2 u &=& 0 &\rm{in}\ \Omega,\\
\frac{\partial u}{\partial \nu}&=&T(k)u&\rm{on}\ \Gamma_R.
\end{array}
\right.
\end{equation}

The variational formulation for \eqref{SchOmega} is to find $k \in \mathbb C$ and $u\in H^1(\Omega)$ such that
\begin{equation}\label{DtN-NEP}
(\nabla u, \nabla v) + (V(x)u, v) - k^2(u,v) - \langle T(k)u,v\rangle_{\Gamma_R} = 0
\end{equation}
for all $v\in H^1(\Omega)$, where $(\cdot,\cdot)$ denotes the $L^2(\Omega)$ inner product and $\langle \cdot,\cdot\rangle$ the $H^{-1/2}(\Gamma_R)$--$H^{1/2}(\Gamma_R)$ duality.

Let $(\cdot,\cdot)_1$ be the $H^1(\Omega)$ inner product. Define the operator $B(k) : H^1(\Omega)\rightarrow H^1(\Omega)$ such that
{\small 
\[
\left(B(k)u,v\right)_1 = (\nabla u, \nabla v) + (V(x)u, v) - k^2(u,v) - \langle T(k)u,v\rangle_{\Gamma_R}
\]
}
for all $v\in H^1(\Omega)$.

The Schr\"{o}dinger resonances are the poles of $B(\cdot)^{-1}$. There is a one-to-one correspondence between resonances and the eigenvalues of $B(\cdot)$ (see Theorem C.10 of \cite{DyatlovZworski2019}). Consequently, we compute the eigenvalue problem for $B(\cdot)$, i.e., find $(k, u\ne 0)\in\mathbb C\times H^1(\Omega)$ such that
\begin{equation}\label{DtN-operator-fun-EVP}
B(k)u=0 \quad {\rm in}\ H^1(\Omega).
\end{equation}

\section{4. Finite Element Discretization}\label{FEM}
Now we turn to the finite element discretization of \eqref{DtN-operator-fun-EVP}. Let $\mathcal{T}_n:=\mathcal{T}_{h_n}$ be a regular triangular mesh for $\Omega$ with mesh size $h_n$. Let $V_n\subseteq H^1(\Omega)$ be the Lagrange finite element space on $\mathcal{T}_n$. Define the discrete operator $B^N_n(k) : V_n \rightarrow V_n$ such that
\begin{eqnarray}
    &&(B^N_n(k)u_n,v_n)_1 = (\nabla u_n, \nabla v_n) +(V(x) u_n, v_n)\qquad  \nonumber \\
 \label{BnN}   &&\qquad \qquad - k^2(u_n,v_n) - \langle T^N(k)u_n,v_n\rangle_{\Gamma_R}
\end{eqnarray}
for all $v_n\in V_n$, where  $T^N(k): H^{1/2}(\Gamma_R)\rightarrow H^{-1/2}(\Gamma_R)$ is the truncated DtN mapping given by

\[
T^N(k)\varphi = \sum_{n=0}^{N}{\!}^{{}^{\,\scriptstyle{\prime}}}\frac{k}{\pi}\frac{{H^{(1)}_n}^{\prime}(kR)}{H^{(1)}_n(kR)}\int_0^{2\pi}\varphi(\phi)\cos(n(\theta-\phi))\text{d}\phi.
\]
The integer $N$ is the truncation order. Similar to \eqref{Tk}, $T^N(k)$ is bounded, i.e.,
\begin{equation}\label{TNk}
    \Vert T^N(k)\varphi\Vert_{H^{-1/2}(\Gamma_R)}\leq C\Vert \varphi\Vert_{H^{1/2}(\Gamma_R)}.
\end{equation}

Let $v_n^i, i=1, 2, \ldots, J$ be the basis function for $V_n$. To obtain the matrix form of \eqref{BnN}, let $S$ be the stiffness matrix with $S^{ij}= (\nabla v^j_n, \nabla v^i_n)$, $M_V$ be the potential matrix with $M_V^{ij}=(V(x)v^j_n, v^i_n)$, $M$ be the mass matrix with $M^{ij}=(v^j_n, v^i_n)$, and $E(k)$ be the matrix due to the truncated DtN mapping with $E^{ij}(k)=\langle T^N(k)v^j_n,v^i_n\rangle_{\Gamma_R}$. The nonlinear algebraic eigenvalue problem corresponding to \eqref{BnN} is to find $k \in \mathbb C$ and a nontrivial vector ${\bf u}$ such that
\begin{equation}\label{AEP}
F(k) {\bf u} = {\bf 0}, 
\end{equation}
where
\[
F(k):= S+M_V-k^2 M - E(k).
\]


Let $\Theta \subset \mathbb C$ be a bounded and connected region such that $\partial \Theta$ is a simple closed curve. Assume that $F(\cdot)$ is holomorphic on $\Theta$ and $F(\cdot)^{-1}$ is meromorphic on $\Theta$ with poles being the eigenvalues of $F(\cdot)$. We present a parallel multi-step spectral indicator method to compute all the eigenvalues in $\Theta$. 

Assuming $F$ has no eigenvalues on $\partial \Theta$, define an operator $P \in \mathbb C^{J, J}$ by
\begin{equation}\label{P}
P=\dfrac{1}{2\pi i}\int_{\partial \Theta}F(z)^{-1}dz,
\end{equation}
which is a projection from $\mathbb C^J$ to the generalized eigenspace associated with all the eigenvalues of $F(\cdot)$ in $\Theta$.
If there are no eigenvalues of $F(\cdot)$ in $\Theta$, then $P = 0$, and $P{\boldsymbol f} = {\boldsymbol 0}$ for all
${\boldsymbol f} \in \mathbb C^J$. Otherwise, $P{\boldsymbol f} \ne {\boldsymbol 0}$ with 
probability $1$ for a random vector ${\boldsymbol f}$. Hence $P{\boldsymbol f}$ can be used to decide if $\Theta$ contains eigenvalues or not. 

We compute the projection $P{\boldsymbol f}$ using the trapezoidal quadrature rule
\begin{equation}\label{XLXf}
P{\boldsymbol f} \approx  \dfrac{1}{2 \pi i} \sum_{j=1}^{N_\omega} \omega_j {\boldsymbol x}_j,
\end{equation}
where $\omega_j$'s are quadrature weights and ${\boldsymbol x}_j$'s are the solutions of the linear systems
\begin{equation}\label{linearsys}
F(z_j){\boldsymbol x}_j = {\boldsymbol f}, \quad j = 1, \ldots, N_\omega,
\end{equation}
where $z_j$'s are the quadrature points.
Define an indicator $I_\Theta$ as
\begin{equation}
I_\Theta = \frac{1}{\sqrt{N_\omega}}\left | \frac{P {\boldsymbol f}|_{N_\omega}}{P {\boldsymbol f}|_{N_\omega/2}} \right|,
\end{equation}
where $P {\boldsymbol f}|_{N_\omega/2}$ is the approximation of $P {\boldsymbol f}$ using $N_\omega/2$ quadrature points $\{z_2, z_4, \ldots, z_{N_\omega}\}$. It is expected that $I_\Theta \approx 1$ if there exist eigenvalues in $\Theta$ and $I_\Theta \approx e^{-CN_\omega/2} (\ll 1)$ for some constant $C$ otherwise. 

Let $N_\omega=32$ be the number of the quadrature points, $tol_{ind} = 0.1$ be a threshold value for the indicator, and $tol_{eps}$ be the accuracy of the eigenvalues.  A parallel multi-step spectral indicator method to compute all the eigenvalues in $\Theta$ is as follows.
\begin{itemize}
\item[-] Given a series of disks $\Theta_i$, with center $\theta_i$ and radius $r$, covering $\Theta$.
\item[1.] Compute eigenvalues in $\Theta$ 
    \begin{itemize}
    \item[1.a.] $L = {\rm ceil}(\log_2(r/tol_{eps}))$.
    \item[1.b.] For $l = 1, \ldots, L$, 
    \begin{itemize}
    	\item[-] Compute the indicators of all the regions of level $l$ in parallel.
	\item[-] Uniformly divide the regions for which the indicators are larger than $tol_{ind}$ into smaller regions.
    \end{itemize}
     \item[1.c.] Use the centers of the regions as approximate eigenvalues $k_i$'s.
    \end{itemize}
\item[2.] Compute the eigenvectors associated to the smallest eigenvalue of $F(k_i)$'s in parallel.
\end{itemize}

In Step 1, we cover $\Theta$ with subregions $\{\Theta_i\}_{i=1}^I$ and compute $I_{\Theta_i}$ to determine if $\Theta_i$ contains eigenvalues. If $I_{\Theta_i} > tol_{ind}$, then $\Theta_i$ is subdivided into smaller regions and these regions are saved for the next level of investigation. The procedure continues until the size of the region is smaller than the precision $tol_{eps}$. Then the centers of these small regions are the approximate eigenvalues.

In Step 2, an approximate eigenvalue $k$ is plugged back into $F(\cdot)$ and compute the eigenvector associated to the smallest eigenvalue of $F(k)$. For more details and Matlab codes, we refer the readers to \cite{XiSun2023}.

\section{5. Numerical Examples}\label{NE}
This section presents various examples to demonstrate the effectiveness of the proposed method.  A regular mesh $\mathcal T_{h_1}$ with mesh size $h_1 \approx 0.05$ is generated for $\Omega$. Then we uniformly refine $\mathcal T_{h_1}$ to obtain $\mathcal T_{h_i}, i=2, \ldots, 5$. The linear Lagrange element is used for $V_h$ for all examples to compute the resonances in 
\[
\Theta := \{a+bi\in\mathbb{C}: \,a \in (-4, 4), \,b \in (-4, -0.5)\}.
\]
The relative error and the convergence order are defined, respectively, as
\[
E_j=\frac{|k^{j}-k^{j+1}|}{|k^{j+1}|},\ \ \ j=1,2,3,4, 
\]
and
\[
\log_2\left(\frac{E_{j}}{E_{j+1}}\right),\ \ \ j=1,2,3.
\]

We take $N=20$ terms in the truncated DtN mapping $T^N(k)$. Numerical experiments indicate that the choice does not significantly affect the accuracy of resonances in $\Theta$, which are relatively small in norm. One does want to increase $N$ if more and/or larger resonances are needed, for which finer meshes are necessary as well.

\begin{figure} [ht!]
\begin{center}
\includegraphics[width=0.4\textwidth]{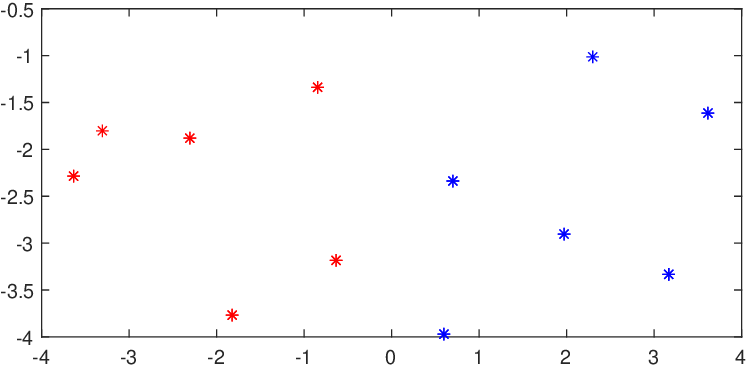}
\caption{Computed resonances of the   ($r_0=1, V(x)=2$).}
\label{figE1}
\end{center}
\end{figure}

{\bf Example 1.}
We first check the case when $D$ is the unit disk and $V(x)=2$ discussed earlier. The computed resonances in $\Theta$ using $\mathcal T_{h_5}$ are shown in Fig.~\ref{figE1}. The values are consistent with Fig.~\ref{PolesSchR1V2}.

In Table~\ref{unitdiskV2}, two resonances (with small absolute values) computed using different meshes are listed. The convergence order is roughly $2$ for both of them, which is expected for the linear Lagrange element.
\begin{table}
\begin{center}
\caption{Two computed resonances and their convergence orders (Example 1).}
\medskip
\label{unitdiskV2}
\begin{tabular}{ccccc}
\hline
 &$k_1$&Ord& $k_2$&Ord\\
\hline
$h_1$   &$-0.835177-1.334674i$&       &$0.588093- 2.144649i$&      \\
$h_2$   &$-0.843659-1.337073i$&       &$0.653423-2.278433i$& \\
$h_3$   &$-0.845785-1.337550i$&2.01 &$0.685785-2.322488i$&1.47 \\
$h_4$   &$-0.846329-1.337659i$&1.97 &$0.695996-2.334152i$&1.82 \\
$h_5$   &$-0.846466-1.337685i$&1.99 &$0.698717-2.337097i$&1.95 \\
\hline
\end{tabular}
\end{center}
\end{table}

\begin{figure}[htb!]
\includegraphics[width=0.4\textwidth]{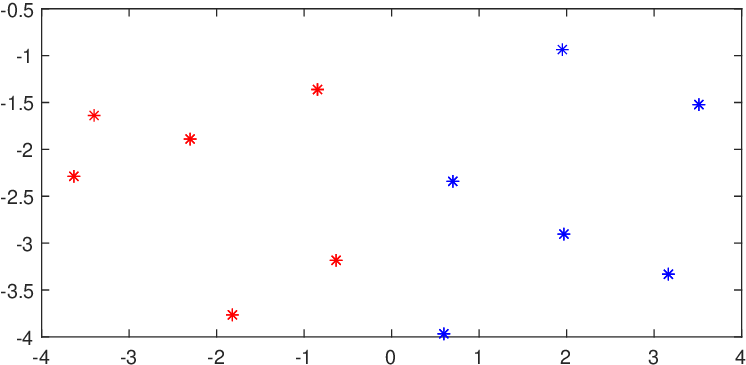} \\
\includegraphics[width=0.4\textwidth]{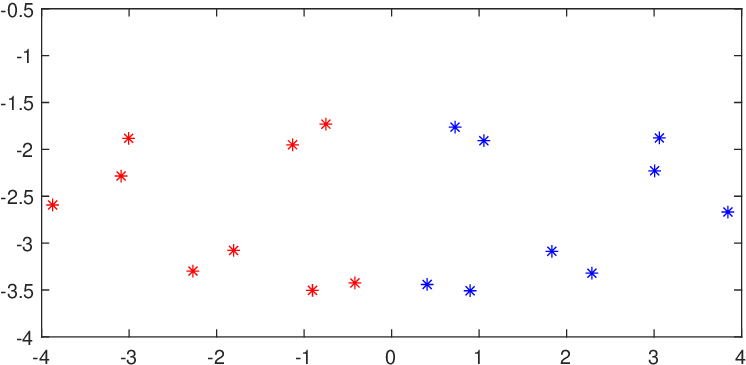} \\
\includegraphics[width=0.4\textwidth]{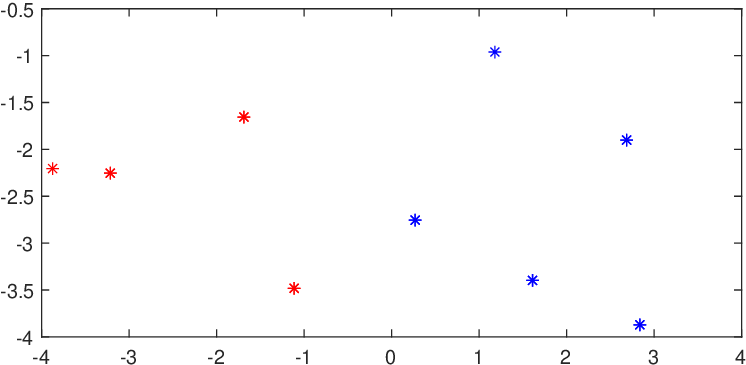} \\
\includegraphics[width=0.4\textwidth]{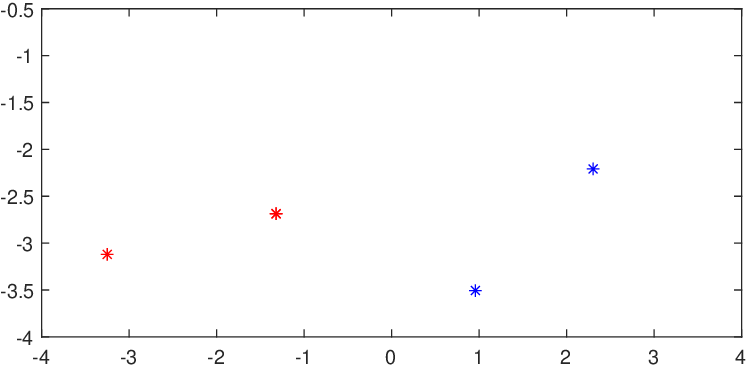}
\caption{\label{resonances} Computed resonances in $\Theta$. First: Example 2. Second: Example 3. Third: Example 4. Fourth: Example 5.}
\end{figure}

{\bf Example 2.}
Let $D:=\{x \in \mathbb R^2\, |\, 1/2 < |x| < 1\}$ and $V(x) = 2$ on $D$.
In Table~\ref{TireV2}, we show two computed resonances and their convergence orders. A second order convergence is obtained. The resonances in $\Theta$ using $\mathcal T_{h_5}$ are shown in Fig.~\ref{resonances} (first picture). 


\begin{table}
\begin{center}
\caption{Two computed resonances and their convergence orders (Example 2).}
\medskip
\label{TireV2}
\begin{tabular}{ccccc}
\hline
 &$k_1$&Ord& $k_2$&Ord\\
\hline
$h_1$   &$-0.835233-1.357152i$&       &$0.589230- 2.147664i$&           \\
$h_2$   &$-0.844288-1.360395i$&       &$0.653609-2.281471i$&     \\
$h_3$   &$-0.846672-1.361085i$&1.96 &$0.685518-2.325891i$&1.47 \\
$h_4$   &$-0.847280-1.361248i$&1.98 &$0.695591-2.337740i$&1.82 \\
$h_5$   &$-0.847433-1.361288i$&1.99 &$0.698266-2.340681i$&1.97 \\
\hline
\end{tabular}
\end{center}
\end{table}

{\bf Example 3.}
Let $D$ be the unit disk and
{\small
\[
V(x) = \left \{ \begin{array}{ll} 2, & x=(r\cos \theta, r\sin \theta), 0 \le r \le 1, \theta \in (0, \pi/2), \\ -1, & x=(r\cos \theta, r\sin \theta), 0 \le r \le 1, \theta \in (\pi/2, \pi), \\1, & x=(r\cos \theta, r\sin \theta), 0 \le r \le 1, \theta \in (\pi, 3\pi/2),\\ -2, & x=(r\cos \theta, r\sin \theta), 0 \le r \le 1, \theta \in (3\pi/2, 2\pi).\end{array}\right.
\]
}
In Table~\ref{Example3}, we show two computed resonances and their convergence orders. The second order convergence is obtained. The resonances in $\Theta$ using $\mathcal T_{h_5}$ are shown in Fig.~\ref{resonances} (second picture). 

\begin{table}
\begin{center}
\caption{Two computed resonances and their convergence orders (Example 3).}
\medskip
\label{Example3}
\begin{tabular}{ccccc}
\hline
 &$k_1$&Ord& $k_2$&Ord\\
\hline
$h_1$   &$-0.734978-1.722162i$&       &$0.701145- 1.799165i$&     \\
$h_2$   &$-0.748316-1.729371i$&       &$0.721810- 1.773786i$&  \\
$h_3$   &$-0.750893-1.730783i$&2.37 &$0.723549-1.765834i$&2.00 \\
$h_4$   &$-0.751487-1.731123i$&2.10 &$0.723781-1.763807i$&1.99 \\
$h_5$   &$-0.751611-1.731209i$&2.18 &$0.723820-1.763287i$&1.98 \\
\hline
\end{tabular}
\end{center}
\end{table}

{\bf Example 4.}
Let $D$ be the unit disk and
\[
V(x) = \exp \left(\frac{1}{|x|^2 -2}\right) + i \exp \left(\frac{1}{|x|^2-4} \right), \quad |x| < 1.
\]
In Table~\ref{Example4}, we show two computed resonances and their convergence orders. The second order convergence is obtained. The resonances in $\Theta$ using $\mathcal T_{h_5}$ are shown in Fig.~\ref{resonances} (third picture). 

\begin{table}[ht!]
\begin{center}
\caption{Two computed resonances and their convergence orders (Example 4).}
\medskip
\label{Example4}
\begin{tabular}{ccccc}
\hline
 &$k_1$&Ord& $k_2$&Ord\\
\hline
$h_1$   &$1.175540- 0.975004i$&       &$-1.772287- 1.632259i$&   \\
$h_2$   &$1.179448-0.964531i$&       &$-1.709829- 1.652785i$&  \\
$h_3$   &$1.179538-0.962084i$&2.19 &$-1.693616-1.654859i$&2.00 \\
$h_4$   &$1.179767-0.961407i$&1.78 &$-1.689543-1.655194i$&2.00 \\
$h_5$   &$1.179824-0.961238i$&2.00 &$-1.688522-1.655267i$&2.00 \\
\hline
\end{tabular}
\end{center}
\end{table}

{\bf Example 5.}
Let $D$ be the union of four squares $D_1=(0.1,0.6) \times (0.1,0.6)$, $D_2= (-0.6,-0.1) \times (0.1,0.6) $, $D_3=(-0.6,-0.1) \times (-0.6,-0.1)$ and $D_4=(0.1,0.6)\times (-0.6,-0.1)$. The potential function $V(x)=2$ on $D$. In Table~\ref{Example5}, we show two computed resonances and their convergence orders. The second order convergence is obtained. The resonances  in $\Theta$ using $\mathcal T_{h_5}$ are shown in Fig.~\ref{resonances} (fourth picture). 

\begin{table}[ht!]
\begin{center}
\caption{Two computed resonances and their convergence orders (Example 5).}
\medskip
\label{Example5}
\begin{tabular}{ccccc}
\hline
 &$k_1$&Ord& $k_2$&Ord\\
\hline
$h_1$   &$-1.222912-2.719221i$&       &$0.768960- 3.213702i$&    \\
$h_2$   &$-1.300488-2.695611i$& &$0.880949-3.423462i$& \\
$h_3$   &$-1.317844-2.685624i$&2.02 &$0.947351-3.494427i$&1.32 \\
$h_4$   &$-1.321996-2.682930i$&2.02 &$0.969666-3.512228i$&1.78\\
$h_5$   &$-1.323022-2.682245i$&2.00 &$0.975692-3.516580i$&1.94\\
\hline
\end{tabular}
\end{center}
\end{table}

\end{document}